\documentclass[12pt]{article}
\usepackage{amsfonts,amsthm,amsmath,amssymb}
\usepackage{graphicx}
\usepackage{wrapfig}

\newtheorem{thm}{Theorem}[section]
\newtheorem{prop}[thm]{Proposition}
\newtheorem{lemma}[thm]{Lemma}
\newtheorem{cor}[thm]{Corollary}
\newtheorem*{example}{Example}

\def\C{\mathbb{C}}
\def\R{\mathbb{R}}
\def\Z{\mathbb{Z}}
\def\P{\mathbb{P}}

\def\Di{\mathbb{D}}

\def\g{\mathfrak{g}}
\def\h{\mathfrak{h}}
\def\p{\mathfrak{p}}

\def\Ocl{\overline{\cal O}}
\def\DG{G\times G}

\def\d{\partial}

\def\t{\widetilde}

\def\a{\alpha}

\def\De{\Delta}

\def\phi{\varphi}

\def\Ga{\Gamma}
\def\L{\Lambda}
\def\la{\lambda}
\def\om{\omega}

\def\Oc{{\cal O}}
\def\F{{\cal F}}

\def\D{{\cal D}}

\def\dim{{\rm dim}}

\def\rk{{\rm rk}}

\def\G{{\rm G}}

\def\End{{\rm End}}

\def\vol{{\rm Vol}}

\def\pic{{\rm Pic}}

\renewcommand\o{\overline}

\title{On intersection indices of subvarieties in reductive groups}
\author{Valentina Kiritchenko}
\date{}

\begin{document}
\maketitle

{\small I will present an explicit formula for the intersection indices of the Chern classes (defined in \cite{VK})
of an arbitrary reductive group with hypersurfaces. This formula has the following applications. First,
it allows to compute explicitly the Euler characteristic of complete intersections in reductive groups.
Second, for any regular compactification of a reductive group, it computes the intersection indices of the Chern
classes of the compactification with hypersurfaces.
The formula is similar to the Brion--Kazarnovskii formula
for the intersection indices of hypersurfaces in reductive groups. The proof uses an algorithm of De Concini and Procesi
for computing such intersection indices. In particular, it is shown that this algorithm produces the Brion--Kazarnovskii
formula.}
\bigskip

\section{Introduction}

Let $G$ be a
connected complex  reductive group of dimension $n$, and let
$\pi:G\to GL(V)$ be a faithful representation of $G$. A generic {\em hyperplane section} $H_\pi$ corresponding
to $\pi$ is the preimage $\pi^{-1}(H)$ of the intersection of $\pi(G)$ with a generic affine hyperplane
$H\subset\End(V)$. There is a nice explicit formula for the self-intersection index of $H_\pi$ in $G$, and more generally,
for the intersection index of $n$ generic hyperplane sections corresponding to different representations
(see Theorem \ref{t.degree} below) in terms of the weight polytopes of the representations.
In this paper, I give a similar formula for the intersection indices of the {\em Chern classes}
of $G$ (defined in \cite{VK}) with generic hyperplane sections (see Theorem \ref{t.chern}).

The {\em Chern classes} of $G$ can  be defined as the Chern classes of the
{\em logarithmic} tangent bundle over a {\em regular}
compactification of $G$ (see Section \ref{s.chern} for a precise definition).
Denote by $k$ the rank of $G$, i.e. the dimension of a maximal torus in $G$. Only the first $(n-k)$
Chern classes are not trivial \cite{VK}. These Chern classes are elements
of the {\em ring of conditions} of $G$, which was introduced by C.De Concini and C.Procesi (see \cite{CP2}).
They can be represented by subvarieties $S_1,\ldots, S_{n-k}\subset G$, where $S_i$ has codimension $i$. All enumerative
problems for $G$, such as the computation of the intersection index $S_iH_\pi^{n-i}$, make sense in the ring
of conditions.

First, I recall the usual Brion--Kazarnovskii formula for the
intersection indices of hyperplane sections. Choose a maximal torus $T\subset G$,
and denote by $L_T$ its character lattice. Choose also a Weyl chamber $\D\subset
L_T\otimes\R$.
Denote by $R^+$ the set of all positive roots of $G$ and denote by $\rho$
the half of the sum of all positive roots of $G$. The inner product
$(\cdot,\cdot)$ on $L_T\otimes\R$ is given by a nondegenerate
 symmetric bilinear form on the Lie algebra of $G$ that is invariant under the adjoint action of $G$  (such a form exists
since $G$ is reductive).
\begin{thm} {\em \cite{Brion, Kaz}} \label{t.degree} If $H_\pi$ is a
hyperplane section corresponding to a
representation $\pi$ with the
weight polytope $P_\pi\subset L_T\otimes\R$ ,
then the self-intersection index of $H_\pi$ in the ring of conditions is equal to
$$n!\int\limits_{P_\pi\cap\D}\prod_{\a\in R^+}\frac{(x,\a)^2}{(\rho,\a)^2}dx.$$
The  measure $dx$ on $L_T\otimes\R$ is normalized so that the covolume of $L_T$ is $1$.
\end{thm}
This theorem was first proved by B.Kazarnovskii \cite{Kaz}. Later, M.Brion proved
an analogous formula for arbitrary spherical varieties using a different method \cite{Brion}.

The integrand in this formula has the following interpretation.
The direct sum $L_T\oplus L_T$ can be identified with the Picard
group of the product $G/B\times G/B$ of two flag varieties. Here
$B$ is a Borel subgroup of $G$. Hence, to each lattice point
$(\la_1,\la_2)\in L_T\oplus L_T$ one can assign the
self-intersection index of the corresponding divisor in $G/B\times
G/B$. The resulting function extends to the polynomial
function $(n-k)!F$ on $(L_T\oplus L_T)\otimes\R$, where
$$F(x,y)=\prod_{\a\in R^+}\frac{(x,\a)(y,\a)}{(\rho,\a)^2}.$$
Note that the integrand is
the restriction of $F$ onto the diagonal $\{(\la,\la):\la\in L_T\otimes\R\}$.

This interpretation leads to another proof of the
Brion--Kazarnovskii formula (different from those of Kazarnovskii and Brion). Namely, take any regular
compactification $X$ of $G$ that lies over the compactification
$X_\pi$ corresponding to the representation $\pi$ (see Subsection \ref{ss.comp}
). Then reduce
the computation of $H_\pi^n$ to the computation of the
intersection indices of divisors in the closed orbits of $X$ (see Section \ref{s.proof}). All
closed orbits are isomorphic to the product of two flag varieties.
The precise algorithm for doing this was given by De Concini and
Procesi \cite{CP} in the case, where $X$ is a wonderful compactification of
a symmetric space. Then E.Bifet extended this algorithm to all regular compactifications of symmetric
spaces \cite{Bifet}. I will show that in the case, where a symmetric space is a reductive group, this
algorithm actually produces the Brion--Kazarnovskii formula if one uses the weight
polytope of $\pi$ to keep track of all transformations.

Moreover, the De Concini--Procesi algorithm
works not only for divisors. It can also be carried over to
the Chern classes
of $G$ (which are, in general, not linear combinations of complete intersections). In particular,
there is the following explicit formula for the intersection
indices of the Chern classes of $G$ with hyperplane sections. Assign to
each lattice point $(\la_1,\la_2)\in L_T\oplus L_T$ the
intersection index of the $i$-th Chern class of the tangent bundle over $G/B\times G/B$
with the divisor $D(\la_1,\la_2)$ corresponding to
$(\la_1,\la_2)$, that is the number $c_i(G/B\times
G/B)D^{n-k-i}(\la_1,\la_2)$. Extend this function to the polynomial
function on $(L_T\oplus L_T)\otimes\R$.
Since the Chern classes of $G/B$ are known the resulting function can be easily computed (see Section \ref{s.proof}).
The final formula is as follows.

Let $\Di$ be the
differential operator (on functions on $(L_T\oplus L_T)\otimes\R$) given by the formula
$$\Di=\prod_{\a\in R^+}(1+\d_\a)(1+\t\d_\a),$$
where $\d_\a$ and $\t\d_\a$ are directional  derivatives along the
vectors $(\a,0)$ and $(0,\a)$, respectively. Denote by $[\Di]_i$
the $i$-th degree term in $\Di$.

\begin{thm} \label{t.chern} If $H_\pi$ is a generic hyperplane section corresponding to a
representation $\pi$ with the weight polytope $P_\pi\subset
L_T\otimes\R$, then the intersection index of $H_\pi^{n-i}$ with
the $i$-th Chern class of $G$ in the ring of conditions is equal
to
$$(n-i)!\int\limits_{P_\pi\cap\D}[\Di]_iF(x,x)dx. \eqno (1)$$
The  measure $dx$ on $L_T\otimes\R$ is normalized so that the covolume of $L_T$ is $1$.
\end{thm}

 Since in general, the Chern classes of $G$ are not complete
intersections, this extends computation of the
intersection indices to a bigger part of the ring of conditions of $G$.
Theorem \ref{t.chern} also completes some results of \cite{VK}. Namely, the Chern classes $S_1,\ldots, S_{n-k}$
were used there as the main ingredients in an adjunction formula for the topological Euler characteristic of
complete intersections of hyperplane sections in $G$ (see Theorem 1.1 in \cite{VK}). Theorem \ref{t.chern}
in the present paper allows to make this formula explicit. E.g. if a complete intersection is just one hyperplane
section $H_\pi$, then
$$\chi(H_\pi)=
(-1)^{n-1}\int\limits_{P_\pi\cap\D}\left(n!-(n-1)![\Di]_1+(n-2)![\Di_2]-\ldots+k![\Di]_{n-k}\right)F(x,x)dx.$$
There is also a formula for the Chern classes $c_i(X)$
of the tangent bundle over any regular compactification  $X$ of $G$  in terms of $S_1$,
\ldots, $S_{n-k}$ (see Corollary 4.4 in \cite{VK}).
Theorem \ref{t.chern} allows to compute explicitly the intersection index of $c_i(X)$ with
a complete intersection of complementary dimension in $X$.

I am grateful to M.Brion, K.Kaveh, B.Kazarnovskii and A.Khovanskii for useful discussions.

\section{Preliminaries}
In this section, I recall some well-known facts which are used in the proof of Theorem \ref{t.chern}.
In Subsection \ref{ss.comp}, I define the
regular compactification $X$ of $G$ associated with a representation $\pi$  and describe the
orbit structure of $X$ in terms of
the weight polytope of the representation. In Subsection \ref{ss.picard}, I will relate the Picard group of $X$
to the space of virtual polytopes
analogous to the weight polytope of $\pi$. The notion of analogous polytopes is discussed in Subsection
\ref{ss.polytopes}.
Then I recall a formula for the integral of a polynomial function
over a simplex (Subsection \ref{ss.integral}), which is used to interpret the computation of intersection
indices in terms of integrals over the weight polytope.

\subsection{Polytopes}\label{ss.polytopes}
Let $P\subset\R^k$ be a convex polytope.  Define the {\em normal fan}
$P^*$ of $P$. This is a fan in the dual space $(\R^k)^*$. To each face $F^i\subset P$ of dimension $i$
there corresponds a cone $F^*_i$ of dimension $(n-i)$ in $P^*$ defined as follows. The cone $F^*_i$
consists of all linear functionals in $(\R^k)^*$ whose maximum value on $P$ is attained on the interior of
the face $F^i$.
In particular, to each facet of $P$ there corresponds a one-dimensional cone, i.e. a ray, in $P^*$.
If the dual space $(\R^k)^*$ is identified with $\R^k$ by means of the Euclidean inner product, the ray
corresponding to a facet is spanned by a normal vector to the facet.

Two convex polytopes are called {\em analogous} if they have the same normal fan. All  polytopes analogous to a given
polytope $P$ form
a semigroup $S_P$ with respect to Minkowski sum. This semigroup is also endowed with the action of the multiplicative
group $\R^{>0}$ (polytopes can be dilated). Hence, $S_P$ can be regarded as a cone in the vector space $V_P$,  where
$V_P$ is the minimal group containing $S_P$ (i.e. the Grothendieck group of $S_P$). The elements of
$V_P$ are called {\em virtual} polytopes analogous to $P$.

We now introduce special coordinates in the vector space $V_P$.
Let $\Ga_1$, \ldots, $\Ga_l$ be the facets of $P$, and let $\Ga_1^*$, \ldots, $\Ga_l^*$ be the corresponding
rays in $P^*$. Choose a non-zero functional $h_i\in\Ga_i^*$ in each ray. Call $h_i$ a {\em support function}
corresponding to the facet $\Ga_i$.
For any polytope $Q$ analogous to $P$, denote by $h_i(Q)$ the maximal value
of $h_i$ on the polytope $Q$. For instance, if $h_i$ is normalized so that its value on the external unit normal to the
facet $\Ga_i$ is 1, then $h_i(P)$ is up to a sign the distance from the origin to the hyperplane that contains the facet
$\Ga_i$ (the sign
is positive if the origin and the polytope $P$ are to the same side of this hyperplane, and negative otherwise). The numbers
$h_1(Q),\ldots,h_l(Q)$ are called the {\em support numbers} of $Q$.
Clearly, the polytope $Q$ is uniquely defined by its support numbers. The coordinates
$h_1(Q),\ldots,h_l(Q)$ can be extended to the space $V_P$,
providing the isomorphism between $V_P$ and the coordinate space $\R^l$.

In what follows, we will deal with {\em integer} polytopes, i.e. polytopes whose vertices belong to a given lattice
$\Z^k\subset \R^k$. For such polytopes, the natural way to normalize the support functions is to require that $h_i(P)$
be equal to the {\em integral} distance from the origin to the hyperplane that contains the facet $\Ga_i$. Suppose that
a hyperplane $H$ not passing through the origin is spanned by lattice vectors. Then the
{\em integral} distance from the origin to the hyperplane $H$  is the index in $\Z^k$ of the subgroup spanned by
$H\cap\Z^k$. To compute the integral distance one can apply a unimodular (with respect to the lattice $\Z^k$)
linear transformation of $\R^{k}$ so
that $H$ becomes parallel to a coordinate hyperplane. Then the integral distance is the usual Euclidean distance from
the origin to this coordinate hyperplane.

We will also use the notion of {\em simple} polytopes.  A polytope in $\R^k$ is called {\em simple} if it is generic
with respect to parallel translations of its facets. Namely, exactly $k$ facets must meet at each vertex.
This implies that any other face is also the transverse intersection of those facets that contain it.

\subsection{Regular compactifications of reductive groups}\label{ss.comp}
With any representation $\pi:G\to GL(V)$ one can associate
the following compactification of $\pi(G)$. Take the projectivization $\P(\pi(G))$ of $\pi(G)$ (i.e. the set
of all lines in $\End(V)$ passing through a point of $\pi(G)$ and the origin),
and then take its closure in $\P(\End(V))$.
We obtain a projective variety $X_\pi\subset\P(\End(V))$
with a natural action of $G\times G$ coming from the left and right action of
$\pi(G)\times\pi(G)$ on $\End(V)$. E.g. when $G=(\C^*)^n$ is a complex torus,
all projective toric varieties can be constructed in this way.

Assume that $\P(\pi(G))$ is isomorphic to $G$.
Consider all weights of the representation $\pi$, i.e. all characters
of the maximal torus $T$ occurring in $\pi$. Take their convex hull $P_\pi$ in $L_T\otimes\R$.
Then it is easy to see that $P_\pi$ is a polytope
invariant under the action of the Weyl group of $G$. It is called the {\em weight polytope}
of the representation $\pi$.
The polytope $P_\pi$ contains
information about the compactification $X_\pi$.
\begin{thm}\label{t.equiv}
1)  {\em (\cite{Tim}, Proposition 8)} The subvariety $X_\pi$ consists of a finite number of
$\DG$-orbits. These orbits are in one-to-one correspondence with
the orbits of the Weyl group acting on the faces of the polytope
$P_\pi$. This correspondence preserves incidence relations. I.e.
if $F_1$, $F_2$ are faces such that $F_1\subset F_2$, then the orbit corresponding to $F_1$ is contained
in the closure of the orbit corresponding to $F_2$.

2)  Let $\sigma$ be another representation of $G$.
The normalizations of subvarieties $X_\pi$ and $X_\sigma$ are isomorphic if and only if
the normal fans corresponding to the polytopes $X_\pi$ and $X_\sigma$
coincide. If the first fan is a subdivision of the second, then there
exists a $\DG$--equivariant map from the normalization of $X_\pi$ to $X_\sigma$, and vice versa.
\end{thm}

The second part of Theorem \ref{t.equiv} follows from the general theory of spherical varieties
(see \cite{Knop}, Theorem 5.1) combined with the description of compactifications $X_\pi$
via colored fans (see \cite{Tim}, Sections 7, 8).

In what follows, we will only consider {\em regular} compactifications of $G$. The simplest example of a regular
compactification is the {\em wonderful compactification} constructed by De Concini and Procesi.
Suppose that the group $G$ is of adjoint type, i.e. the center of $G$ is trivial.
Take  any irreducible representation $\pi$  with a strictly dominant highest weight.
It is proved in \cite{CP} that the corresponding compactification $X_\pi$ of the group $G$
is always smooth and, hence, does not depend on the choice of a highest
weight. Indeed, the normal fan of
the weight polytope $P_\pi$ coincides with the fan of the Weyl chambers and their faces, so the second
part of Theorem \ref{t.equiv}
applies.
This compactification is
called the {\em wonderful compactification} and is denoted by $X_{can}$.

Other {\em regular} compactifications of $G$ can be characterized as follows.
The normalization $X$ of $X_\pi$ is {\em regular} if first, it is smooth, and second, there is
a $(\DG)$--equivariant map from $X$ to $X_{can}$.
These two conditions can be reformulated in terms of the weight polytope $P_\pi$.
Namely, the first condition implies that $P_\pi$ is {\em integrally} simple (see \cite{Tim} Theorem 9), i.e. it is simple
and the edges meeting at each vertex form a basis of $L_T$. The second condition
implies that none of the vertices of $P_\pi$ lies on the
walls of the Weyl chambers, i.e. the normal fan of $P_\pi$ subdivides the fan of the Weyl chambers
and their faces.

A regular compactification $X$ has the following nice properties (see \cite{Brion2} for details), which we will use
in the sequel.
The  boundary divisor $X\setminus G$ is a divisor with normal crossings.
The $\DG$--orbits of codimension $s$ correspond to the faces of $P_\pi$ of codimension $s$ and have rank $(k-s)$.
Recall that each face $F\subset P_\pi$ is the transverse intersection of several facets of $P_\pi$ (since $P_\pi$
is simple). Then the closure of the orbit corresponding to $F$  is the transverse intersection of the closures
of the codimension one orbits that correspond to these facets. Each closed orbit of $X$ (such orbits correspond to the
vertices of $P_\pi$) is isomorphic to the product of
two flag varieties $G/B\times G/B$.

\subsection{Picard groups of compactifications}\label{ss.picard}
Let $X$ be the normalization of the compactification $X_\pi$ of $G$.
We assume that $X$ is regular, and hence smooth. Then the second cohomology group $H^2(X)$ is isomorphic to the Picard group of $X$
(see \cite{Bifet}).
 Denote by $V(\pi)$
the group of all integer
virtual polytopes analogous to the weight polytope $P_\pi$ and invariant under the action of the Weyl group.
There is a description of the Picard group of a regular complete symmetric space due to Bifet
(see \cite{Bifet}, Theorem 2.4, see also \cite{Brion}, Proposition 3.2).
In our case, this description can be reformulated as follows (such a reformulation is well-known in the toric case, and
in the reductive case it was suggested by K.Kiumars).
The Picard group $\pic(X)$ of $X$ is canonically isomorphic to the quotient group of $V(\pi)$ modulo
parallel translations.
In particular, if $G$ is semisimple, then $\pic(X)=V(\pi)$ (the only parallel translation
taking a $W$--invariant polytope to a $W$--invariant polytope is the trivial one).
The isomorphism takes the hyperplane section corresponding to a representation $\sigma$ to the weight polytope
of $\sigma$ and extends to the other divisors by linearity. Let us identify divisors in $X$ with the corresponding
polytopes using this isomorphism.

The variety $X$ has $l$ distinguished {\em boundary divisors} $\Ocl_1$, \ldots, $\Ocl_l$, which are the closures of
codimension one orbits. Let us describe the corresponding virtual polytopes.
Choose $l$ facets $\Ga_1$, \ldots,
$\Ga_l$ of $P_\pi$ so that each orbit of the Weyl group acting  on the facets of $P_\pi$
contains exactly one $\Ga_i$. E.g. take all facets that intersect the fundamental Weyl chamber.
Choose the support functions $h_1$, \ldots, $h_l$ corresponding to these facets
so that $h_i(P_\pi)$ is equal to the integral distance (with respect to the weight lattice $L_T$) from
the origin to the facet $\Ga_i$.

\begin{lemma} \label{l.orbits} The closure $\Ocl_i$ of codimension one orbit corresponds to the virtual polytope
whose $i$-th support number is $1$ and the other support numbers are $0$.
\end{lemma}
\begin{proof}
Let $\sigma$ be any representation of $G$ whose weight polytope $P$ is analogous to $P_\pi$. Then $X$ is isomorphic to
the normalization of
the compactification $X_\sigma$. Thus a generic linear functional $f$ on $X_\sigma$ can
also be regarded as a rational function on $X$. Let us find the zero and the pole divisors of $f$.
The zero divisor is the divisor corresponding to the weight polytope of $\sigma$. The pole divisor
is a linear
combination of the divisors $\Ocl_1$, \ldots, $\Ocl_l$. It is not hard to show that the coefficients are
the support numbers $h_1(P_\sigma)$, \ldots, $h_l(P_\sigma)$, i.e. the integral distances from the origin to
the facets of $P_\sigma$
corresponding to $\Ga_1$, \ldots, $\Ga_l$. Indeed, for toric varieties, this statement is well-known
(see \cite{Fulton}, Section 3.4). In particular, this holds
for the closure $\o T$ in $X$ of the maximal torus $T\subset G$. Hence,
this also holds for $X$, and the divisor $D$ of a hyperplane section corresponding to $\sigma$
can be written as
$$D=h_1(P_\sigma)\Ocl_1+\ldots+h_l(P_\sigma)\Ocl_l.$$
It follows that $h_i(\Ocl_j)=0$, unless $i=j$.
\end{proof}
Another useful collection  of divisors  consists of the closures in $X$
of codimension one Bruhat cells in $G$.  Denote these divisors by $D_1$, \ldots, $D_k$.
They can also be described as hyperplane sections corresponding to the irreducible representations of $G$
with  fundamental highest weights $\om_1$, \ldots, $\om_k$, respectively.
Then to each dominant
weight $\la=m_1\om_1+\ldots+m_k\om_k$ there corresponds the {\em weight divisor} $D(\la)=m_1D_1+\ldots+m_kD_k$.
The polytope of this divisor is the weight polytope $P_{\la}$ of the irreducible representation with the
highest weight $\la$.  Note that $\la$ is the only vertex of $P_{\la}$ inside the fundamental Weyl chamber.
Hence, it belongs to all facets of $P_{\la}$ corresponding to $\Ga_1$,\ldots, $\Ga_l$
(e.g. some of the facets might degenerate to the vertex $\la$). This implies the following lemma.

\begin{lemma} \label{l.cells} Let $D(\la)$ be the weight divisor corresponding to a
weight $\la\in L_T$ and let $P_\la$ be its polytope.
Then $h_i(P_\la)=h_i(\la)$ for any $i=1,\ldots,l$.
\end{lemma}

Combination of these two lemmas leads to the following result.

\begin{cor}\label{c.divisor}
Let $D$ be the divisor on $X$ corresponding to a polytope $P$. We assume that $P$ is analogous to $P_\pi$ and identify
the respective facets.
Then for any face
$F\subset P$ of codimension $s$ that intersects the fundamental Weyl chamber $\D$
 and for any point $\la\in F\cap\D$, the divisor $D$ can be written uniquely as a linear combination of $D(\la)$
and of boundary divisors $\Ocl_i$ such that the corresponding facets $\Ga_i$ do not contain $F$. Namely, if
$F=\Ga_{i_1}\cap\ldots\cap\Ga_{i_s}$, then
$$D=D(\la)+\sum_{j\in\{1,\ldots,l\}\setminus\{i_1,\ldots,i_s\}}[h_j(P)-h_j(\la)]\Ocl_j.$$
\end{cor}

\subsection{Integration of polynomials} \label{ss.integral}

Let $f(x_1,\ldots,x_k)$ be a homogeneous polynomial function of
degree $d$ defined on a real affine space $\R^k$ with coordinates
$(x_1,\ldots,x_k)$. Below I recall a useful formula expressing the
integral of $f$ over a simplex in $\R^k$ in terms of the
{\em polarization} of $f$. Recall that the {\em polarization} of $f$ is the unique
symmetric $d$-linear form $f_{pol}$ on $\R^k$ such that the restriction of $f_{pol}$
to the diagonal coincides with $f$. One can define $f_{pol}$ explicitly as follows:
$$f_{pol}(v_1,\ldots,v_d)=\frac{1}{d!}\frac{\d^d}{\d_{v_1}\ldots\d_{v_d}}f,$$
where $\d_{v_i}$ is the directional derivative along the vector $v_i$.

 Let $\Delta\subset\R^k$ be a $k$-dimensional simplex with
vertices $a_0,\ldots,a_k$ and let $dx=dx_1\wedge
dx_2\wedge\ldots\wedge dx_k$  be the standard measure on $\R^k$.

\begin{prop}{\em \cite{Brion}}\label{p.int}
Let $f_{pol}$ be the polarization of $f$. It can be regarded as a
linear function on  the $d$-th symmetric power of $V$. Then the
average value of $f$ on the simplex $\Delta$ coincides with the
average value of $f_{pol}$ on all symmetric products of $d$
vectors from the set \{$a_0$, \ldots, $a_k$\}:
$$\frac1{\vol(\Delta)}{\int\limits_\Delta f(x)dx}=
\frac1{\binom{d+k}{k}}
\sum_{i_0+\ldots+i_k=d}
f_{pol}(\underbrace{a_0,\ldots,a_0}_{i_0},\ldots,\underbrace{a_k,\ldots,a_k}_{i_k}).$$
\end{prop}

\section{Chern classes}\label{s.chern}
In this section, I recall the definition of the Chern classes of spherical homogeneous spaces (see \cite{VK}
for more details). In the sequel, only Chern classes of $\DG$--orbits in regular compactifications of $G$
will be used.
For these Chern classes, I prove a vanishing result for their intersection indices with
certain weight divisors in regular compactifications. This result will be important in Section \ref{s.proof} when
applying the De Concini--Procesi algorithm to the Chern classes of $G$.

Let $G/H$ be a spherical homogeneous space under $G$. Denote by $\g$ and $\h$ the Lie algebras of $G$ and $H$,
respectively, and denote by $m$ the dimension of $\h$.
Define the {\em Demazure map} $p$ from $G/H$ to the Grassmannian $\G(m,\g)$ of $m$-dimensional subspaces in $\g$
as follows:
$$p:G/H\to\G(m,\g); \quad p: gH\to g\h g^{-1}.$$
Let $C_i\subset\G(m,\g)$ be the Schubert cycle corresponding to a generic subspace $\L_i\subset\g$
of codimension $m+i-1$, i.e. $C_i=\{\L\in\G(m,\g): \dim(\L\cap\L_i)\ge1\}$. Then the $i$-th {\em  Chern class}
$S_i(G/H)$ of $G/H$ is the preimage of $C_i$ under the map $p$:
$$S_i(G/H)=p^{-1}(C_i).$$
The class of $S_i(G/H)$ in the ring of conditions of $G/H$ is the same for all generic $C_i$ \cite{VK}.
It is related to the Chern classes of the tangent bundles over  regular compactifications of $G$
\cite{Brion3,VK}. Namely, if $X$ is a regular compactification of $G/H$, then the closure of $S_i(G/H)$
in $X$ is the $i$-th Chern class of the {\em logarithmic} tangent bundle over $X$ that corresponds
to the divisor $X\setminus (G/H)$. This vector bundle is generated by all vector fields on $X$ that are tangent
to $G$--orbits in $X$. In what follows, this bundle will be called the {\em Demazure bundle} of $X$.

Let $X=G/H$ and $Y=G/P$ be two spherical homogeneous
spaces under $G$.  Suppose that $H$ is a subgroup of $P$. Consider the $G$--equivariant map
$$p:X\to Y; \quad p: gH\mapsto gP.$$
Let $S_i(X)$ the $i$-th Chern class of $X$.
In general, it is not true that $S_i(X)$ is the inverse image under the map $p$
of a subset in $Y$. However, the intersection of $S_i(X)$ (when it is nonempty) with a fiber of $p$ has
dimension at least $\rk(P)-\rk(H)$.

{\begin{example}\em In what follows, we will mostly deal with the case, where $X$ and $Y$ are spherical
homogeneous spaces
under the doubled group $\DG$. Namely, $X$ is a $\DG$--orbit of a regular compactification
of the group $G$
and $Y$ is a partial flag variety constructed as
follows.  Let $H\subset \DG$ be the stabilizer of a
point in $\Oc$. Take the minimal parabolic subgroup $P\subset G\times G$ that contains $H$ and set $F=(\DG)/P$.
\end{example}}
\begin{lemma}\label{l.fibers}
For a generic $S_i(X)$, there exists an open dense subset of $S_i(X)$  such that for any element $x$ of this subset
the intersection of the fiber $xP$ with $S_i(X)$ has dimension
greater than or equal to the $\rk(P)-\rk(H)$. In particular, the dimension of $p(S_i(X))$ satisfies
the inequality
$$\dim~p(S_i(X))\le\dim~S_i(X)-(\rk(P)-\rk(H)).$$
\end{lemma}

\begin{proof}
Choose a generic vector space  $\L\subset\g$ of codimension $\dim~H+i-1$.
Denote by $\h$ and $\p$ the Lie algebras of $H$ and $P$
respectively. Then by definition
$S_i(X)$ consists of all cosets $gH$ such that $g\h g^{-1}$ has a nontrivial intersection with $\L$,
or equivalently $\h\cap g^{-1}\L g$ is nontrivial.

Let $gH$ be any element of $S_i(X)$. Estimate the dimension of the intersection of $S_i(X)$ with
the fiber $gP$ of the map $p$.
Note that for all $g$ from a dense open subset of $S_i(X)$, the intersection
$\h\cap g^{-1}\L g$ contains an element $v$ that is regular in  $\h$.
Denote by $C$ the centralizer
in $P$ of $v\in\h\subset\p$. Then $\dim(C\cap H)=\rk(H)$ while $C$ has dimension at least $\rk(P)$.
Note that for any $c\in C$ the coset $gcH$ still belongs to
$S_i(X)$ since $c^{-1}g^{-1}\L gc$ contains $c^{-1}vc=v$.
Hence, $S_i(X)\cap gP$ contains a set $gCH$ of dimension
at least $\rk(P)-\rk(H)$.
\end{proof}
Lemma \ref{l.fibers} is crucial for proving the following two vanishing results,
which extend Proposition 9.1 from \cite{CP} and rely on the same ideas.
Let $X$ be a regular compactification of $G$, and let $p:X\to X_{can}$ be its equivariant projection to the wonderful
compactification. Denote by $c_1$, \ldots, $c_{n-k}$ the Chern classes of the Demazure vector bundle
over $X$.

\begin{lemma} \label{l.vanish}
Let $\Oc$ be a $\DG$--orbit in $X$ of codimension $s<k$, and $\Ocl\subset X$ its closure.
Suppose that the image $p(\Oc)$ under the map $p:X\to X_{can}$
coincides with the closed orbit of $X_{can}$. In terms of polytopes, this means that the face corresponding to
$\Oc$  does not intersect the walls
of the Weyl chambers.

Let $\la$ be any weight of $G$, and $D(\la)$ the corresponding weight divisor. Then the homology class
$c_iD^{n-i-s}(\la)$ vanishes on $\Ocl$, i.e. the following intersection index is zero:
$$c_iD(\la)^{n-i-s}\Ocl=0.$$
\end{lemma}

\begin{proof} First of all, the intersection product $c_i\cdot\Ocl$ is the $i$-th Chern class
of the Demazure bundle over $\Ocl$ (see \cite{Bien}, Proposition 2.4.2). Hence, it can be realized as the
closure in $\Ocl$ of the
$i$-th Chern class $S_i(\Oc)$ of the spherical homogeneous space $\Oc$.
The computation of the intersection index
$c_i\Ocl D(\la)^{n-i-s}$ in $X$ thus reduces to the computation of the intersection index
$\o{S_i(\Oc)}D(\la)^{n-i-s}$ in $\Ocl$. The latter is equal to the intersection index
$S_i(\Oc)D(\la)^{n-i-s}$ in the ring of conditions of $\Oc$ since $D(\la)$ and $\o{S_i(\Oc)}$
have proper intersections with the boundary $\Ocl\setminus \Oc$.

To compute $S_i(\Oc)D(\la)^{n-i-s}$ we use the restriction of the map $p:X\to X_{can}$
to $\Ocl$. By the hypothesis the image $p(\Ocl)$ is the closed orbit $F$ in $X_{can}$, so it is isomorphic to
the product $G/B\times G/B$ of two flag varieties. Then the divisor $D(\la)$
restricted to $\Ocl$ is the inverse image under the map $p$ of the divisor $D(\la,\la)$ in $F$. Indeed,
$D(\la)=p^{-1}(\t D(\la))$, where $\t D(\la)$ is the weight divisor in $X_{can}$ corresponding to $\la$.
It is easy to check that $\t D(\la)\cap F=D(\la,\la)$ (see Proposition 8.1 in \cite{CP}).

Hence, all the intersection points in
$S_i(\Oc)D(\la)^{n-i-s}$ are contained in the preimage of
$p(S_i(\Oc)) D(\la,\la)^{n-i-s}$. But the latter is empty.
Indeed, since $\Oc$ has positive rank and $F$ has zero rank,
Lemma \ref{l.fibers} implies that
$$\dim~p(S_i(\Oc))<\dim~S_i(\Oc)=n-i-s.$$

\end{proof}

It remains to deal with the orbits in $X$ whose image under the map $p$ is not the closed orbit in $X_{can}$.
In this case, the face corresponding to such an orbit intersects the walls of the Weyl chambers, and hence, it
is orthogonal to some of the fundamental weights $\om_1$,\ldots, $\om_k$.  Note that the
codimension one orbits $\Oc_1$,\ldots, $\Oc_k$ in $X_{can}$ are in one-to-one correspondence with the
fundamental weights $\om_1$, \ldots, $\om_k$. Namely, the facet corresponding to $\Oc_i$ is orthogonal to $\om_i$.
Let $\Ocl_1$,\ldots, $\Ocl_k$ be the closures in $X_{can}$ of $\Oc_1$,\ldots, $\Oc_k$, respectively.
\begin{lemma}
Let $\Oc$ be a $\DG$--orbit in $X$ of codimension $s<k$.
Suppose that the image $p(\Oc)$ under the map $p:X\to X_{can}$
is not closed and lies in the intersection $\Ocl_{i_1}\cap\ldots\cap \Ocl_{i_s}$. In terms of polytopes, this means
that the face corresponding to $\Oc$ is  orthogonal to the weights $\om_{i_1}$,\ldots, $\om_{i_s}$.

Let $\la$ be any linear combination of the weights $\om_{i_1}$,\ldots, $\om_{i_s}$.  Then
$$c_iD(\la)^{n-i-s}\Ocl=0.$$
\end{lemma}
\begin{proof}
We use the $\DG$--equivariant map $r$ from $\Ocl_{i_1}\cap\ldots\cap \Ocl_{i_s}$
to a partial flag variety $G/P\times G/P$ constructed in \cite{CP} (see \cite{CP} Lemma 5.1 for details).
Consider the compactification $X_{i_1,\ldots,i_s}$
of $G$ corresponding to the irreducible representation $\pi_{i_1,\ldots,i_s}$ whose highest weight lies
strictly inside the cone spanned by $\om_{i_1}$,\ldots, $\om_{i_s}$. This compactification has a unique closed
orbit $G/P\times G/P$, where $P\subset G$ is the stabilizer of the highest weight vector in the representation
$\pi_{i_1,\ldots,i_s}$.
Clearly, the fan of the Weyl chambers and their faces subdivides the normal fan of the weight polytope of
$\pi_{i_1,\ldots,i_s}$. Hence, by Theorem \ref{t.equiv} there is an equivariant map $r:X_{can}\to X_{i_1,\ldots,i_s}$.
This map takes $\Ocl_{i_1}\cap\ldots\cap \Ocl_{i_s}$ to the closed orbit $G/P\times G/P$.

The composition $rp$ maps the orbit $\Oc$ to the closed orbit $G/P\times G/P$ of $X_{i_1,\ldots,i_s}$.
It is easy to show that
the divisor $D(\la)$ restricted to $\Oc$ is the preimage of the divisor $D(\la,\la)\subset G/P\times G/P$ under this map
(see \cite{CP} Section 8.1).
Now repeat the arguments of the proof of Lemma \ref{l.vanish}.

\end{proof}

These two lemmas imply the following vanishing result.

\begin{cor}\label{c.vanish}
Let $\Oc$ be any $\DG$--orbit in $X$, and let $F$ be the face of the polytope of $X$ that corresponds to $\Oc$.
The intersection index
$$c_i D(\la_1)\ldots D(\la_{n-i-s})\Ocl$$
vanishes in the cohomology ring of $X$ in the following two cases:

1) The face $F$ does not intersect the walls of the Weyl chambers. Then weights $\la_1$,\ldots, $\la_{n-i-s}$
are any weights of $G$.

2) The face $F$ intersects a wall of the Weyl chambers and weights $\la_1$,\ldots, $\la_{n-i-s}$ are
orthogonal to $F$ (with respect to the inner product $(\cdot,\cdot)$ on $L_T\otimes\R$ defined in the Introduction).
\end{cor}

\section{Proof of Theorem \ref{t.chern}} \label{s.proof}
We use notation of Subsections \ref{ss.comp} and \ref{ss.picard}.
Let $X$ be any
regular compactification lying over the compactification $X_\pi$. Then the closure $\o H_\pi$ of $H_\pi$ in $X$
has proper intersections with all $\DG$--orbits in $X$, and thus $S_iH^{n-i}_\pi$
coincides with the intersection index $\o S_i\o H_\pi^{n-i}$ in the cohomology ring of $X$.

Assume that $X$ corresponds to a representation of $G$ with the weight polytope
$P_0$. Let us compute $\o S_iD^{n-i}$ for a divisor $D$ under the  assumption that the polytope $P$
corresponding to $D$ is analogous to $P_0$. After we establish the formula of Theorem \ref{t.chern} for such divisors,
it will automatically extend to the other divisors (in particular, for $\o H_\pi$) since any virtual polytope analogous
to $P_0$ is a linear combination of polytopes analogous to $P_0$. Since $X$ is regular, $P_0$ and hence $P$ are simple.

All computations are carried in the cohomology ring of $X$. First, break
$D^{n-i}$ into monomials of the form
$\Ocl_{i_1}\ldots\Ocl_{i_k}D(\la_1)\ldots D(\la_{n-i-k})$, where
$i_1$,\ldots, $i_k$ are distinct integers from $1$ to $l$ and
$\la_1$,\ldots, $\la_{n-i-k}$ are  weights. Then every such monomial can be computed explicitly,
since the intersection $\Ocl_{i_1}\cap\ldots\cap\Ocl_{i_k}$ is either empty or isomorphic to the
product of two flag varieties.

\begin{wrapfigure}{l}{6.5cm}
\includegraphics[scale=.4]{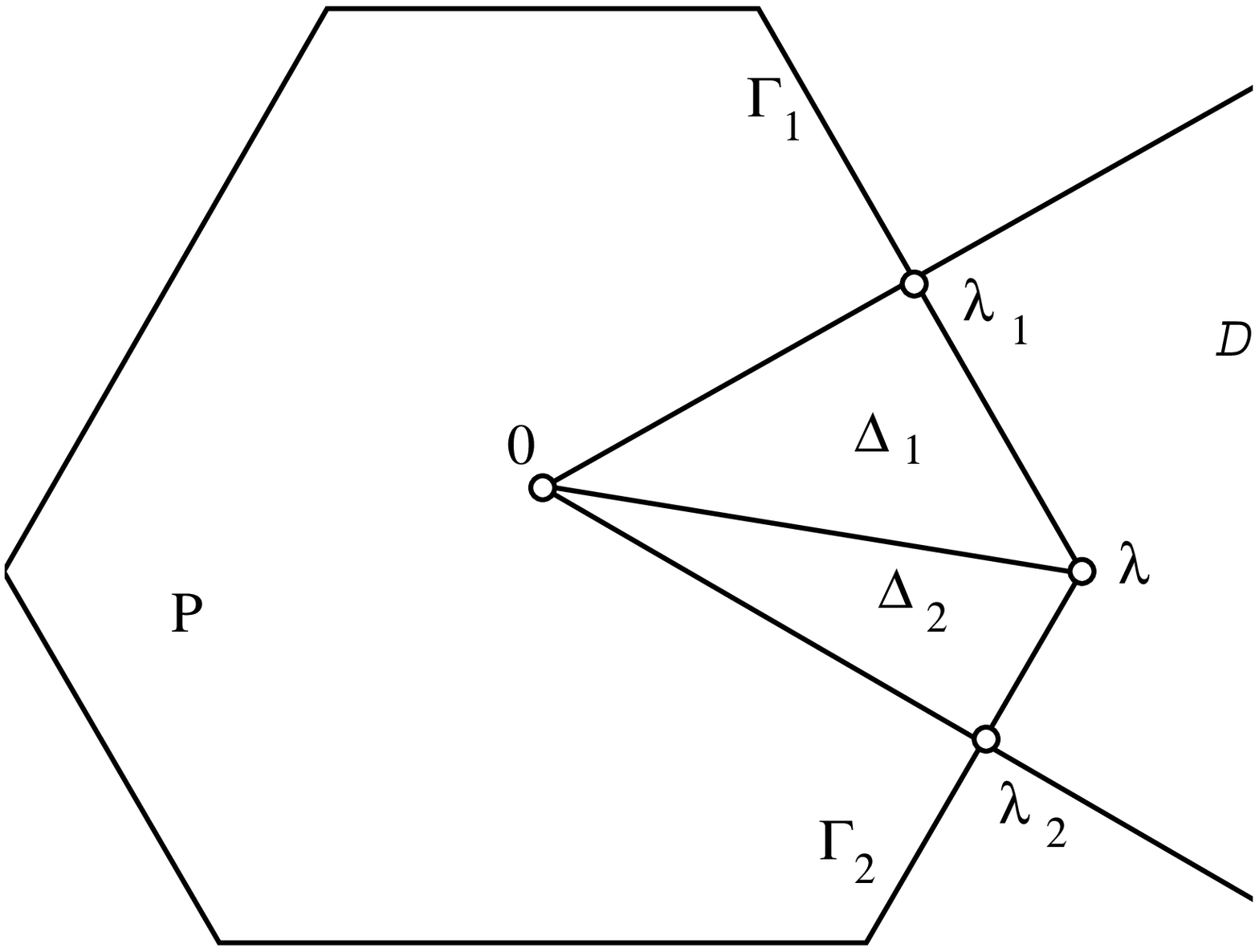}
\caption{}
\end{wrapfigure}

Since we are going to intersect $D^{n-i}$ with $\o S_i$
we can ignore all monomials that are annihilated by $\o S_i$. Recall that $\o S_i$ is the $i$-th Chern class of the
Demazure bundle over $X$. In particular,
Corollary \ref{c.vanish} implies that $\o S_i$ annihilates the ideal $I\subset H^*(X)$
generated by the monomials of the form $D(\la_1)\ldots D(\la_{n-i-s})\Ocl$ such that
either the face of $P$ corresponding to the codimension $s$ orbit $\Oc$ does not intersect the walls of the
Weyl chambers or, if it does, the weights $\la_1,\ldots,\la_{n-i-s}$ are orthogonal to this face.

To keep track of our calculations we use a subdivision of the
polytope $P\cap\D$ into simplices coming from the barycentric subdivision of $P$ described below.
For each face $F\subset P$  choose a point $\la_F\in F$ as follows. If $F$ does not intersect the
walls of the Weyl chamber $\D$, then $\la_F$ is any  point in the interior of  the face. Otherwise, choose $\la_F$
so that the corresponding vector is orthogonal to the face $F$ (in particular, $\la_F$ will belong to
the intersection of the face with a wall of $\D$.)
If $F=P$ take $\la_F=0$.

An {\em $s$--flag} $\F$ is the collection $\{F_1\supset\ldots\supset F_s\}$ of $s\le k$
nested faces of $P$ such that each of them intersects
$\D$, and $F_i$ has codimension $i$ in $P$.
Denote by $\Ocl_\F$ the closure in $X$ of the orbit corresponding
to the last face $F_s$, and by $\De_\F$ the $s$--dimensional simplex with the vertices $0$, $\la_{F_1}$,
\ldots, $\la_{F_s}$. In particular, when $s=k$, the simplex $\De_\F$ has full dimension and the orbit $\Ocl_\F$
is closed. The polytope $\D\cap P$ is the union of simplices $\De_\F$ over all
possible $k$--flags $\F$.

{\bf Example.} Take $G=PSL_3(\C)$, and let $X=X_{can}$ be
its wonderful compactification. Let divisor $D$ be a hyperplane section corresponding
to the irreducible representation with a strictly dominant highest weight $\lambda$.
In this case, $P$ is a hexagon symmetric under the action of the Weyl group with
two edges $\Ga_1$ and $\Ga_2$ intersecting $\D$.
Then $\la_i=\la_{\Ga_i}\in\Ga_i$ is orthogonal to $\Ga_i$ for $i=1,2$ and $\Ga_1\cap\Ga_2=\la$.
The subdivision of $P\cap\D$ into simplices
consists of two triangles $\De_1$ and $\De_2$ with the vertices $0$, $\la_1$, $\la$ and $0$, $\la_2$, $\la$,
respectively (see Figure 1).

\begin{lemma} \label{l.equiv} Denote by $f_d(x_1,\ldots,x_k)$  the sum of all
monomials of degree $d$ in $k$ variables $x_1$,\ldots, $x_k$.
The following identity holds in the cohomology ring of $X$ modulo the ideal $I$:
$$D^{n-i}\equiv k!\sum_{\F}\vol(\De_\F)
f_{n-k-i}(D,D(\la_{F_1}),\ldots,D(\la_{F_{k-1}}))\Ocl_{\F}\pmod I,$$
where the sum is taken over all possible $k$--flags $\F=\{F_1\supset\ldots\supset F_k\}$.
The volume form $\vol$ is normalized so that the covolume of $L_T$ is equal to $1$.
\end{lemma}

\begin{proof} We will prove the following more general statement for $s$-flags. Denote by
$f_{d,s}(x_1,\ldots,x_s)$ the sum of all monomials of degree $d$ in $s$ variables.

Recall that $\Ga_1$,\ldots, $\Ga_l$ denote the facets of $P$
that intersect the Weyl chamber $\D$. An $s$-flag  can be alternatively described by an ordered collection of
facets $\Ga_{i_1}$,\ldots, $\Ga_{i_s}$ such that their intersection $\Ga_{i_1}\cap\ldots\cap\Ga_{i_s}$ has codimension
$s$. Then $F_j=\Ga_{i_1}\cap\ldots\cap\Ga_{i_j}$. This is a one-to-one correspondence, since the polytope $P$ is simple.
Assign to each $s$-flag $\F$ the following number
$$c_\F=h_{i_1}(P)[h_{i_2}(P)-h_{i_2}(\la_{F_1})]\ldots[h_{i_s}(P)-h_{i_s}(\la_{F_s})].$$

In particular, when $s=k$, i.e. $F_s$ is just a vertex, the number $c_\F$
coincides with the volume of $\De_\F$ times $k!$. Indeed, by a unimodular linear transformation of $L_T\otimes\R$
we can map the hyperplanes containing the facets $\Ga_{i_1}$,\ldots, $\Ga_{i_s}$ to the coordinate hyperplanes.
Then $[h_{i_j}(P)-h_{i_j}(\la_{F_{j-1}})]$ is just the Euclidean distance from the vertex $\la_{F_{j-1}}$ of $\De_\F$
to the hyperplane containing $\Ga_{i_j}$. Note that to define volumes we do not use the inner product $(\cdot,\cdot)$ on
the lattice $L_T$ defined in the introduction. We only use the lattice itself.

Then for any integer $s$ such that $1\le s\le k$ the following is true:
$$D^{n-i}\equiv\sum_{\F}c_\F
f_{n-s-i,s}(D,D(\la_{F_1}),\ldots,D(\la_{F_{s-1}}))\Ocl_{\F}\pmod I, \eqno(1)$$
where the sum is taken over all $s$--flags.

Prove by induction on $s$. We use the notations of Subsection \ref{ss.picard}.
For $s=1$, the statement coincides with the decomposition
$D=h_1(P)\Ocl_1+\cdots+h_l(P)\Ocl_l$ from Lemma \ref{l.orbits}.

Assume that the formula is proved for some $s<k$. Prove it for $s+1$. We now deal separately with each term on
the right hand side of formula (1). First subtract from every term
$f_{n-s-i,s}(D,D(\la_{F_1}),\ldots,D(\la_{F_{s-1}}))\Ocl_\F $
the element $f_{n-s-i,s}(D(\la_{F_s}),D(\la_{\F_1}),\ldots,D(\la_{F_{s-1}}))\Ocl_\F$ of the ideal $I$.
This operation does not change the identity (1). A simple calculation shows that
$$f_{n-s-i,s}(x,x_1,\ldots,x_{s-1})-f_{n-s-i,s}(x_s,x_1,\ldots,x_{s-1})
=(x-x_s)f_{n-s-i-1,s+1}(x,x_1,\ldots,x_{s-1},x_s).$$
Hence, after subtraction we can rewrite the difference as
$$(D-D(\la_{F_s}))f_{n-s-i-1,s+1}(D,D(\la_{F_1}),\ldots,D(\la_{F_s}))\Ocl_\F.$$
Since $\la_s$ lies in the intersection of $s$ facets
 $\Ga_{i_1}$,\ldots, $\Ga_{i_s}$, Corollary \ref{c.divisor} implies that
$$(D-D(\la_{F_s}))\Ocl_\F=\sum_{j\ne i_1,\ldots,i_k}[h_j(P)-h_j(\la_{F_s})]\Ocl_j\Ocl_\F.$$
Note that $\Ocl_j\Ocl_\F$ is empty if and only if the intersection of $\Ga_j$ with $\Ga_{i_1}\cap\ldots\cap\Ga_{i_s}$ is
empty.
Hence,
$$(D-D(\la_{F_s}))\Ocl_\F=\sum_{\F'}[h_j(P)-h_j(\la_{F_s})]\Ocl_{\F'},$$
where the sum is taken over all $(s+1)$-flags $\F'$ that extend $\F$, i.e.
$\F'=\{F_1\supset\ldots\supset F_s\supset F_s\cap\Ga_{j}\}$.
\end{proof}

It remains to compute the term
$$\o S_i\cdot f_{n-k-i}(D,D(\la_{F_1}),\ldots,D(\la_{F_{k-1}}))\Ocl_{\F} \eqno(2)$$
for each $k$-flag $\F$. Suppose that the closed orbit $\Ocl_\F$ is the intersection of $k$ hypersurfaces
$\Ocl_{i_1}$,\ldots, $\Ocl_{i_k}$. Then for any other codimension 1 orbit $\Oc_j$ (such that $j\ne i_1,\ldots,i_k$),
the intersection $\Ocl_\F\cap\Ocl_j$ is empty. Hence, $D$ in (2) can be replaced by
$D(\la_{F_k})$ since
$$D=D(\la_{F_k})+\sum_{j\ne i_1,\ldots,i_k}(h_j(P)-h_j(\la_{F_k}))\Ocl_j.$$
Note also that the evaluation of (2) reduces to the computation of
intersection indices in $\Ocl_\F$, which is the product of two flag varieties. We have that
$\o S_i\cdot\Ocl_\F=c_i(\Ocl_\F)$ and $D(\la)\cdot\Ocl_\F=D(\la,\la)$. Here $c_i(\Ocl_\F)$ is the $i$-th Chern class
of the tangent bundle over $\Ocl_\F$. Hence,
$$\o S_if_{n-k-i}(D(\la_{F_k}),D(\la_{F_1}),\ldots,D(\la_{F_{k-1}}))\Ocl_{\F}=$$
$$=c_i(G/B\times G/B)f_{n-k-i}(D(\la_{F_1},\la_{F_1}),\ldots,D(\la_{F_k},\la_{F_k})). \eqno (3)$$
The intersection product in the right hand side of this formula is taken in $G/B\times G/B$.

The function $F_i(\la)=c_i(G/B\times G/B)D(\la,\la)^{n-k-i}$ can be expressed
explicitly in terms
of the function $F$ defined in the Introduction, since the $i$-th Chern class of
$G/B\times G/B$ is the term of degree $i$ in the intersection product
$$\prod_{\a\in R^+}(1+D(\a,0))(1+D(0,\a)).$$
One way to compute $F_i$ is as follows.
Let $\Di$ and $[\Di]_i$ be the
differential operators defined in the Introduction. Then
$$F_i(x)=(n-k-i)![\Di]_iF(x,x).$$
This easily follows from the formula for the polarization mentioned in Subsection \ref{ss.integral}
and the fact that $D^{n-k}(\la,\la)=(n-k)!F(\la,\la)$.

We can now apply Proposition \ref{p.int} to convert the sum (3) into the integral over the simplex $\De_\F$.
Indeed, by definition of the function $f_{n-k-i}$ we have that (3) can be rewritten as
$$\sum_{i_1+\ldots+i_k=n-k-i}
{F_i}_{pol}(\underbrace{\la_{F_1},\ldots,\la_{F_1}}_{i_1},\ldots,\underbrace{\la_{F_k},\ldots,\la_{F_k}}_{i_k}).$$
This is equal to the integral
$$\binom{n-i}{k}{\int\limits_{\Delta_\F} F_i(x)dx}/{\vol(\Delta_\F)}$$
by Proposition \ref{p.int} applied to the simplex $\De_\F$ (with the vertices $0$, $\la_{F_1}$,\ldots, $\la_{F_k}$)
and to the function $F_i(x)$.
Combining this with Lemma \ref{l.equiv} we get
$$\o S_iD^{n-i}=\frac{(n-i)!}{(n-k-i)!}
\sum_{\F}{\int\limits_{\Delta_\F} F_i(x)dx}=(n-i)!\int\limits_{P\cap\D}[\Di]_iF(x,x)dx.$$
Note that when $i=0$, we get the Brion--Kazarnovskii formula.

\end{document}